\documentclass[12pt]{article}

\usepackage{amsmath}
\usepackage{theorem}
\usepackage{amssymb}
\usepackage{latexsym}
\usepackage{a4wide}
\usepackage{longtable}
\usepackage{epic}

\newtheorem{thm}{Theorem}
\newtheorem{lem}[thm]{Lemma}
\newtheorem{prop}[thm]{Proposition}

\newtheorem{conj}[thm]{Conjecture}
{\theorembodyfont{\rmfamily} }
\newenvironment{rem}{\noindent{\bf Remark.}}{\newline}
\newenvironment{pf}{\noindent{\bf Proof.}}{\hbox{}\hfill $\Box$}

\newcommand{\Q}{\mathbb{Q}}
\newcommand{\R}{\mathbb{R}}
\newcommand{\C}{\mathbb{C}}
\newcommand{\F}{\mathbb{F}}
\newcommand{\Z}{\mathbb{Z}}
\newcommand{\Aut}{\mathrm{\mathop{Aut}}}
\newcommand{\Lie}{\mathrm{\mathop{Lie}}}
\newcommand{\Gal}{\mathrm{\mathop{Gal}}}
\newcommand{\GL}{\mathrm{\mathop{GL}}}
\newcommand{\Tr}{\mathrm{\mathop{Tr}}}
\newcommand{\Ad}{\mathrm{\mathop{Ad}}}
\newcommand{\ad}{\mathrm{\mathop{ad}}}
\newcommand{\ind}{\mathrm{\mathop{ind}}}
\newcommand{\rank}{\mathrm{\mathop{rank}}}
\newcommand{\gl}{\mathfrak{\mathop{gl}}}
\newcommand{\mf}[1]{\mathfrak{#1}}
\newcommand{\ssl}{\mathfrak{\mathop{sl}}}
\newcommand{\g}{\mathfrak{g}}
\newcommand{\lab}[1]{\footnotesize{#1}}

\begin{document}

\title{Computing with nilpotent orbits in simple Lie algebras of exceptional type}
\author{Willem A. de Graaf\\
Dipartimento di Matematica\\
Universit\`{a} di Trento\\
Italy}
\date{}
\maketitle

\begin{abstract}
Let $G$ be a simple algebraic group over an algebraically closed field with Lie algebra 
$\g$. Then the orbits of nilpotent elements of $\g$ under the adjoint action of $G$ 
have been classified. We describe a simple algorithm for finding a representative of 
a nilpotent orbit. We use this to compute lists of representatives of these orbits
for the Lie algebras of exceptional type.
Then we give two applications. The first one concerns settling a conjecture by 
Elashvili on the index of centralizers of nilpotent orbits, for the case where the 
Lie algebra is of exceptional type. The second deals with minimal dimensions of
centralizers in centralizers.
\end{abstract}

\section{Introduction}\label{sec:intro}

Let $G$ be a simple algebraic group over an algebraically closed field of characteristic $0$.
Let $\g$ denote its Lie algebra. Then $G$ acts on $\g$ via the adjoint representation. 
It is a natural question what the $G$-orbits in $\g$ are. Recall that an element $e\in \g$
is said to be nilpotent if the map $\ad e : \g\to \g$ is nilpotent. Now the $G$-orbits of
nilpotent elements in $\g$ are called nilpotent orbits. These have drawn a lot of attention
in the past decades. On some occasions it turns out that using conceptual arguments
to prove their properties is a lot harder for the exceptional types than it is for the
classical types. However, for the former an approach based on a case by case analysis
is possible. It is the objective of this paper to describe how this can be carried out
using computer calculations.\par
The nilpotent orbits in $\g$ are classified in terms of so-called weighted Dynkin diagrams. 
The first problem that we consider is to find a nilpotent element in $\g$ 
given the corresponding weighted Dynkin diagram. We describe a straightforward 
algorithm for this (Section \ref{sec:orb}). Then the algorithm is used to compute lists of 
explicit representatives of the nilpotent orbits in the Lie algebras of exceptional
type. They are listed in Appendix \ref{sec:tables}.\par
We use these lists to prove Elashvili's conjecture
for the exceptional types by computer calculations. This conjecture concerns the index
of centralizers of nilpotent elements. The concept of index is defined as follows.
Let $K$ be a finite-dimensional Lie algebra, and let $K^*$ denote the dual space.
For $f\in K^*$ set $K^f = \{ x\in K \mid f([x,y])=0 \text{ for all } y\in K\}$.
Then the index of $K$ is defined as the number
$$\ind(K) = \inf_{f\in K^*} \dim K^f.$$
For semisimple Lie algebras in characteristic zero it is known that the index
is equal to the rank (\cite{dix}, Proposition 1.11.12).\par 
By $C_\g(x)$ we denote the centralizer of of $x\in \g$.

\begin{conj}[Elashvili]
Let $\g$ be a semisimple Lie algebra over an algebraically closed field of characteristic
$0$. Let $x\in \g$. Then $\ind( C_\g(x) )$ is equal to the rank of $\g$.
\end{conj}

This conjecture has recently received renewed attention, cf.
\cite{panpreya}, \cite{panyushev}, \cite{vinya}. Its proof 
immediately reduces to the case where 
$\g$ is simple, and $x$ nilpotent (cf. \cite{panyushev}, \S 3). Also an 
inequality of Vinberg states that
$\ind( C_\g(x) )$ is at least the rank of $\g$ (see \cite{panyushev}, 1.6, 1.7).
The conjecture has been proved for $\g$ of classical type in \cite{yakimova}, see also the
discussion in \cite{panpreya}. In Section \ref{sec:index} we report on computer calculations 
that settle the conjecture for the exceptional types. \par
In \cite{sekiguchi} the question is considered whether for a given nilpotent $e\in \g$ there exists 
$x\in C_\g(e)$ such that the dimension of $C_\g(e,x)$ equals the rank 
of $\g$. There an example is given where such an $x$ does not exist, for the case where $\g$ 
is of type $F_4$. In Section \ref{sec:dblcen} we approach this question using our lists of
representatives of nilpotent orbits. This way we are able to give a complete list of all 
$e$ for which such an $x$ does not exist, in all exceptional types. For the Lie algebra
of type $E_8$ this solves an open problem from \cite{sekiguchi}. For type $G_2$ this corrects a 
statement in \cite{sekiguchi}.\par
The paper ends with two appendices. The first contains the lists of representatives of nilpotent
orbits. The second (Appendix \ref{sec:appB}) has lists of positive roots as they appear in
the computer algebra system {\sf GAP}. They have been added to help reading the tables
of Appendix \ref{sec:tables}. \par
All algorithms described in this paper have been implemented
in the language of the computer algebra system {\sf GAP}4. The implementations
are available from
\begin{center}
\begin{verbatim}
http://www.science.unitn.it/~degraaf/nilpotent-orbit.html
\end{verbatim}
\end{center}
{\bf Acknowledgments:} I thank Alexander Elashvili for suggesting
all the topics of this paper to me, and for his enthusiastic advice while I was 
writing it. Also I would like to thank Karin Baur for several helpful email
exchanges, and for her comments on earlier versions.

\section{Preliminaries on nilpotent orbits}\label{sec:prelim}

In this section we give a short overview of the theory behind the classification
of nilpotent orbits. For more detailed accounts we refer to \cite{cart}, \cite{colmcgov}.\par
Let $e\in\g$ be a nilpotent element. Then by the Jacobson-Morozov theorem
$e$ lies in a subalgebra of $\g$ that is isomorphic to $\ssl_2$. In other words, there are
elements $f,h\in \g$ with $[e,f]=h$, $[h,f]=-2f$, $[h,e]=2e$. In this case we say
that $(f,h,e)$ is an $\ssl_2$-triple. \par
Now let $(f,h,e)$ be an $\ssl_2$-triple. Then by the representation theory of $\ssl_2$
we get a direct sum decomposition $\g = \oplus_{k\in \Z} \g(k)$, where 
$\g(k) = \{ x\in \g \mid [h,x]=kx \}$. Fix a Cartan subalgebra $H$ of $\g$ with 
$h\in H$. Let $\Phi$ be the corresponding root system of $\g$.  
For $\alpha\in \Phi$ we let $x_\alpha$ be a corresponding root vector. For each
$\alpha$ there is a $k\in \Z$ with $x_\alpha\in \g(k)$. We write $\eta(\alpha)=k$.
It can be shown that there exists a basis of simple roots $\Delta\subset \Phi$ such that
$\eta(\alpha)\geq 0$ for all $\alpha\in\Delta$. Furthermore, for such a $\Delta$ we
have $\eta(\alpha)\in \{0,1,2\}$ for all $\alpha\in\Delta$. Write $\Delta = \{\alpha_1,
\ldots,\alpha_l\}$. Then the Dynkin diagram of $\Phi$ has $l$ nodes, the $i$-th node
corresponding to $\alpha_i$. Now to each node we add the label $\eta(\alpha_i)$; the
result is called the weighted Dynkin diagram. It is denoted $\Delta(e)$, and it depends 
only on $e$, and not on the choice of $\ssl_2$-triple containing $e$. \par
Let $e,e'$ be two nilpotent elements in $\g$. It can be shown that $e,e'$ lie in
the same $G$-orbit if and only if $\Delta(e)=\Delta(e')$. So the weighted
Dynkin diagram of $e$ uniquely identifies the nilpotent orbit $Ge$. The 
weighted Dynkin diagrams corresponding to nilpotent orbits have been classified. 
For the exceptional types there are
explicit lists. For the classical types there is a classification in terms of 
partitions. In particular, the nilpotent orbits in $\g$ have been classified.\par
Let $e\in \g$ be a representative of a nilpotent orbit. We may assume that $e$ is
a linear combination of root vectors, corresponding to positive roots. Let
$\beta_1,\ldots,\beta_r$ be the positive roots involved in this linear combination.
Let $x_{\beta_i}$ (respectively $y_{\beta_i}$) be the root vector corresponding to 
$\beta_i$ (respectively $-\beta_i$). Let $\mf{l}\subset \g$ be the subalgebra generated
by $H$ along with the $x_{\beta_i}$ and $y_{\beta_i}$. Then $\mf{l}$ is reductive, and 
$e\in \mf{l}$. Let $(f,h,e)$ be an $\ssl_2$-triple containing $e$, contained in $\mf{l}$.
Then $\mf{l}$ decomposes with respect to the action of $\ad h$ as 
$\mf{l} = \oplus_{k\in \Z} \mf{l}(k)$.  Let $\mf{p}=\oplus_{k\geq 0} \mf{l}(k)$, which is a 
subalgebra of $\mf{l}$. Now it can be shown that the nilpotent orbit containing $e$ is
uniquely determined by the pair $(\mf{l},\mf{p})$ (cf. \cite{colmcgov}, Chapter 8).
Corresponding to this the nilpotent orbit has the label $X_n(a_i)$, where $X_n$ is
the type of the semisimple part of $\mf{l}$, and $i$ is the number of simple roots in
the semisimple part of $\mf{p}$. If the latter algebra is solvable, then we omit the
$a_i$. Furthermore, if the roots of $\mf{l}$ are short (seen as roots of $\g$), then
a tilde is put over the $X_n$. On some occasions, two different orbits can have the same
label. Then a $'$ is added to one of them, whereas the other gets $''$.
We note that, although the pair $(\mf{l},\mf{p})$
uniquely determines the nilpotent orbit, it is also true that the same nilpotent orbit
can have more than one (non-isomorphic) such pair. So the same nilpotent orbit can have
more than one label.\par
The nilpotent element $e$ from above also has a Dynkin diagram, which is simply the Dynkin
diagram of the roots $\beta_i$. This diagram has $r$ nodes, and node $i$ is connected
to node $j$ by $\langle \beta_i, \beta_j^\vee\rangle \langle \beta_j,\beta_i^\vee\rangle =0,1,2,3$
lines. Furthermore, if these scalar products are positive, then the lines are dotted.
This only occurs when $\mf{p}$ is not solvable.

\section{Finding representatives of nilpotent orbits}\label{sec:orb}

In this section we consider the problem of finding a nilpotent element in $\g$
corresponding to a given weighted Dynkin diagram $D$. We write $D_i$ for the
label at node $i$. Let $H$ be a fixed Cartan subalgebra of $\g$.\par
Let $e\in \g$ be a nilpotent element such that $\Delta(e)=D$. Then there is an
$\ssl_2$-triple $(f,h,e)$, containing $e$. Since we can conjugate any Cartan
subalgebra of $\g$ to $H$ by an element of $G$, we may assume that $h\in H$.
As in the previous section we write $\Delta=\{\alpha_1,\ldots,\alpha_l\}$ for
a basis of simple roots. By choosing a Chevalley basis in $\g$ we get basis
elements $h_1,\ldots,h_l$ of $H$, and root vectors $x_{\alpha_i}$ with 
$[h_j,x_{\alpha_i}]=\langle \alpha_i,\alpha_j^\vee\rangle x_{\alpha_i}$. \par
Each $h\in H$ yields a decomposition $\g = \oplus_{k\in \Z} \g(k)$, and a
weighted Dynkin diagram, as described in the previous section. This weighted
Dynkin diagram is equal to $D$ if and only if $[h,x_{\alpha_i}]= D_i x_{\alpha_i}$
for $1\leq i\leq l$. But this happens if and only if $\sum_{j=1}^l \langle 
\alpha_i,\alpha_j^\vee\rangle a_j =D_i$, where the $a_j$ are such that 
$h=\sum_j a_j h_j$. Let $C=(\langle \alpha_i,\alpha_j^\vee \rangle)_{1\leq i,j\leq l}$
be the Cartan matrix of $\Phi$. It follows that $h$ yields the weighted Dynkin
diagram $D$ if and only if $C(a_1,\ldots,a_l)^t = (D_1,\ldots,D_l)$. Hence
that there is a unique such $h$, and we can compute it by solving a system of 
linear equations. However, not every weighted Dynkin diagram corresponds to a
nilpotent orbit. In other words, not every weighted Dynkin diagram yields a $h$
that lies in an $\ssl_2$-triple. The next two lemmas lead to a probabilistic algorithm 
to decide whether this is the case or not.

\begin{lem}
Let $h\in H$. Then $h$ belongs to an $\ssl_2$-triple if and only if there is an
$x\in \g(2)$ such that $h\in [x,\g(-2)]$.
\end{lem}

\begin{pf}
The condition is clearly necessary. If $h\in [x,\g(-2)]$ then there is a
$y\in \g(-2)$ with $[x,y]=h$. Then $(y,h,x)$ is an $\ssl_2$-triple. 
\end{pf}

\begin{lem}\label{lem2}
Let $h\in H$ be contained in an $\ssl_2$-triple $(y,h,x)$. Let $E$ be the
set of $x'\in \g(2)$ such that $h\in [x',\g(-2)]$. Then $E$ is Zariski dense
in $\g(2)$.
\end{lem}

\begin{pf}
(cf. \cite{cart}, Proposition 5.6.2). Let $G_h = \{g\in G\mid \Ad(g)(h)=h\}$ be
the stabilizer of $h$ in $G$. Then $G_h$ is an algebraic subgroup of $G$.
Now $\Lie(G_h) = \{u\in \g \mid \ad(u)(h)=0\}$. This is the centralizer of 
$h$ in $\g$. Hence $\Lie(G_h) = \g(0)$. For $u\in \g(2)$ and $g\in G_h$ we
have $[h,\Ad(g)(u)] = \Ad(g)[\Ad(g^{-1})(h),u] = \Ad(g)[h,u]=2\Ad(g)(u)$.
Hence $\Ad(g)$ stabilizes $\g(2)$. Let $\varphi : G_h \to \g(2)$ be the morphism defined
by $\varphi(g)=\Ad(g)(x)$. Then the image of $\varphi$ is the $G_h$-orbit
of $x$ in $\g(2)$. The differential of $\varphi$ is $d\varphi : \g(0)\to \g(2)$,
$d\varphi(u) = [u,x]$. But this is surjective because $[\g(0),x]=\g(2)$ (this
follows from the representation theory of $\ssl2_2$). 
So $\varphi$ is a dominant morphism. Hence $\varphi(G_h)$ is a dense subset 
of $\g(2)$. Furthermore $\varphi(G_h)\subset E$. 
\end{pf}

Based on this we have a probabilistic algorithm for finding a representative
of a nilpotent orbit, given a weighted Dynkin diagram. First we determine
the unique $h\in H$ corresponding to the diagram. Then we select a random 
$x\in \g(2)$, in the following way. Let $x_1,\ldots,x_s$ be a basis of $\g(2)$.
Let $\Omega$ be a finite subset of 
$\Q$ and select $\mu_1,\ldots,\mu_s$ randomly, uniformly and
independently from $\Omega$. Then set $x= \sum_i \mu_i x_i$.
 By the previous lemma the probability that $h\in [x,\g(-2)]$
is high (and can be made arbitrarily close to $1$ by enlarging $\Omega$). If it happens to be
the case that $h\not\in [x,\g(-2)]$ then we select another $x$ and continue.
This algorithm will terminate in very few steps. \par
The $x$ found by the algorithm above will have ``ugly'' coefficients with 
respect to a Chevalley basis. We can obtain an element with ``nice'' coefficients
in the following way. We write $x$ with respect to a Chevalley basis of $\g$.
We fix every coefficient but the first. For the first coefficient we try the 
values $0,1,2,\ldots$. The lemma ensures that we will quickly find an $x'$ 
which is a representative of the same nilpotent orbit, with the first coefficient
a nice integer. We continue this way until all coefficients are nice integers.\par
The above results also provide a probabilistic algorithm for testing whether 
a given weighted Dynkin diagram corresponds to a nilpotent orbit. We basically try
the same algorithm a few times, and if it does not come up with an $x$ then 
the weighted Dynkin diagram does not correspond to a nilpotent orbit with
high probability. In principle we can make this absolutely sure by using Gr\"obner bases. 
This works as follows.
Let $x_1,\ldots,x_s$ and $y_1,\ldots,y_s$ be bases of respectively $\g(2)$ and $\g(-2)$.
Let $a_1,\ldots,a_s,b_1,\ldots,b_s$ be indeterminates. Let $u_1,\ldots,u_r$ be a basis
of $\g(0)$, and write $[x_i,y_j] = \sum_k \gamma_{ij}^k u_k$, and $h= \sum_k \alpha_k u_k$.
Then there is an $x\in \g(2)$ with $h\in [x,\g(-2)]$ if and only if the system 
of polynomial equations
$$\sum_{i=1}^s \sum_{j=1}^s \gamma_{ij}^k a_ib_j -\alpha_k = 0 \text{ ~~ for } 1\leq k\leq r$$
has a solution. Now this system has a solution over $\C$ if and only if the reduced
Gr\"obner basis of the ideal generated by the left hand sides of these equations is
not $\{1\}$. \par
\begin{rem}
In \cite{popov} Popov has given an algorithm for determining the strata of 
the nullcone of a linear representation of a reductive algebraic group. This 
also yields an algorithm for classifying nilpotent orbits in reductive Lie algebras,
and for finding representatives of them. 
\end{rem}

\section{Calculating the index}\label{sec:index}

In this section we describe a simple algorithm that for a Lie algebra 
gives an upper bound for its index. If the Lie algebra is defined over
a sufficiently large field (e.g., of characteristic $0$), then the 
probability that this upper bound is 
equal to the index can be made arbitrarily high. We use the same notation
as in Section \ref{sec:intro}. \par 
Let $K$ be a finite-dimensional Lie algebra with basis $\{x_1,\ldots,x_n\}$.
Let $c_{ij}^k$ be the structure constants of $K$, i.e., $[x_i,x_j] = \sum_{k=1}^n
c_{ij}^k x_k$.
Let $\{\psi_1,\ldots,\psi_n\}$ be the dual basis of $K^*$, i.e., $\psi_i(x_j)=
\delta_{ij}$. Let $f=\sum_i T_i \psi_i$ be an element of the dual space $K^*$.
Let $x = \sum_i \alpha_i x_i\in K$. Then $x\in K^f$ if and only if $f([x,x_j])=0$
for $1\leq j\leq n$. Now this is equivalent to
$$\sum_{i=1}^n (\sum_{k=1}^n c_{ij}^k T_k)\alpha_i = 0 \text{ for } j=1,\ldots,n.$$
Define the $n\times n$-matrix $A$ by $A(i,j) = \sum_{k=1}^n c_{ij}^k T_k$.
Then $\dim K^f = n -\rank(A)$. So the dimension of $K^f$ is minimal if and only
if the rank of $A$ is maximal. Now the rank of $A$ is not maximal if and only
if certain polynomial expressions in the $T_k$ 
(i.e., determinants of certain minors of $A$) vanish. Therefore, if the $T_k$ are
chosen randomly and uniformly from a sufficiently large set, then with high
probability the rank of $A$ will be maximal. \par
Here we consider the case where $K=C_\g(e)$, where $e$ is a nilpotent element of
the simple Lie algebra $\g$. Then by Vinberg's inequality we have that 
$\ind(K)$ is at least the rank of $\g$. So if we find an $f$ such that 
$\dim (K^f) = \rank(\g)$, then we have proved that $\ind(K)=\rank(\g)$. 
Moreover, the above discussion shows that we will quickly find such an $f$ 
by randomly choosing the $T_k$. \par
With the help of an implementation of this algorithm in {\sf GAP}4,
we have checked Elashvili's conjecture for the exceptional
types (which, except $G_2$, are the remaining open cases). As a result we can conclude that 
Elashvili's conjecture holds for all simple Lie algebras. \par

\section{Centralizers in centralizers}\label{sec:dblcen}

Let $e\in\g$ be a nilpotent element. Let $C_e= C_\g(e)$ be the centralizer of $e$ in $\g$.
Let $x\in C_e$ and consider the centralizer $C_{e,x}$ of $x$ in $C_e$ (i.e., $C_{e,x}$ is
the set of all elements of $\g$ commuting with both $e$ and $x$). 
From \cite{richardson} it follows that $C_{e,x}$ contains a
commutative subalgebra of dimension equal to $\rank(\g)$. Hence the dimension of $C_{e,x}$
is at least the rank of $\g$. In \cite{sekiguchi} the following question is considered:
given $e$ does there exist $x\in C_e$ such that the dimension of $C_{e,x}$ equals the rank 
of $\g$? The main result of that paper is a counter example to the question for the
case where $\g$ is of type $F_4$. \par
With the lists of representatives of the nilpotent orbits we can easily tackle this
question in all simple Lie algebras of exceptional type. Let $e\in\g$ be a nilpotent
element, and let $x_1,\ldots,x_m$ be a basis of $C_e$. Set $x=T_1x_1+\cdots +T_mx_m$.
Then the centralizer of $x$ in $C_e$ is equal to the kernel of $\ad x$ (restricted to $C_e$).
So the dimension of $C_{e,x}$ is minimal if the rank of the matrix $\ad x$ is maximal.
Now the entries of this matrix are linear polynomials in the $T_i$. It follows that for
a random choice of the $T_i$, with very high probability, the rank of $\ad x$ is maximal.
So this gives a probabilistic algorithm for determining the minimal dimension of $C_{e,x}$
(recall that we are varying $x$, and keeping $e$ fixed). Once the minimal dimension
is found with this algorithm we can prove it rigorously as follows. Let $x$ be an
element such that $\dim C_{e,x}$ is (hypothetically) minimal, as produced by the algorithm. 
If $\dim C_{e,x} = \rank(\g)$ then we have proved that the minimal dimension 
of a $C_{e,x}$ is $\rank(\g)$, as it cannot be smaller. 
Secondly, if the dimension that we find happens to be bigger, then we compute
the rank of the matrix $\ad x$, where $x=T_1x_1+\cdots +T_mx_m$ and we let
the $T_i$ be generators of a rational function field. The rank of that matrix will
equal the maximal rank of any $\ad x$ for $x\in C_e$. \par
Using this algorithm we arrive at the following result.

\begin{prop}
Let $\g$ be a simple Lie algebra of exceptional type, and $e\in\g$ nilpotent. Then
the minimal dimension of a $C_{e,x}$ is equal to $\rank(\g)$, except in
three cases, which are listed in the following table:
\begin{longtable}{|l|l|l|}
\endfirsthead 
\hline
\endhead
\hline
\endfoot
\endlastfoot

\hline
type of $\g$ & label of $e$ & dimension of minimal $C_{e,x}$ \\
\hline
$G_2$ & $A_1+\widetilde{A}_1$ & $3$ \\
$F_4$ & $\widetilde{A_2}+A_2$ & $6$ \\
$E_8$ & $A_5+A_2+A_1$ & $12$ \\	
\hline
\end{longtable}
In all three cases it turns out that a minimal $C_{e,x}$ 
is abelian. Furthermore, in each case it is possible to choose the element $x\in C_e$
such that it is homogeneous of degree $-1$ with respect to the grading of $\g$ defined by
the $\ssl_2$-triple containing $e$.
\end{prop}

In relation to \cite{sekiguchi} we remark the following. In \cite{sekiguchi} it is
wrongly stated that in $G_2$ all minimal $C_{e,x}$ have dimension
equal to $\rank(\g)$. The result for $F_4$ is the same as in \cite{sekiguchi}.
Finally, the problem for $E_8$ is left open in \cite{sekiguchi}. \par
Also, as a straightforward corollary of the proposition, it follows that in the
exceptional types a minimal $C_{e,x}$ is always abelian.

\appendix
\section{Representatives of nilpotent orbits}\label{sec:tables}

In the tables below we list the nilpotent orbits in the Lie algebras of exceptional
type. For each orbit we have given a label, the weighted Dynkin diagram, and the
Dynkin diagram of a representative. We remark the following. If more than one label 
was possible, we have chosen the simplest one that we could find. This means that we
have preferred a label of the form $X_n$ over a label of the form $X_n(a_i)$. Furthermore,
we have preferred labels such that the Dynkin diagram of a corresponding 
representative has as few lines as possible. In the Dynkin diagram a black node
means that the corresponding root is long. Finally, the labels corresponding to each
node refer to the basis elements of the simple Lie algebras as present in {\sf GAP}4.
In Appendix \ref{sec:appB} we list the positive roots of each root system of
exceptional type, in the order in which they are used by {\sf GAP}4. Now, if in
the tables in this section a Dynkin diagram of a representative has labels
$i_1,\ldots,i_k$, then the corresponding representative is the sum of the root
vectors corresponding to the $i_j$-th positive root for $1\leq j\leq k$.

\setlongtables
\begin{longtable}{|l|c|l|}
\caption{Nilpotent orbits in the Lie algebra of type $G_2$.}
\endfirsthead
\hline
\multicolumn{3}{|l|}{\small\slshape Nilpotent orbits in type $G_2$} \\
\hline 
\endhead
\hline 
\endfoot
\endlastfoot

\hline
label & diagram & representative \\

      & \begin{picture}(10,7)
  \put(-8,0){\circle*{6}}
  \put(22,0){\circle{6}}
  \put(-7,-3){\line(1,0){28}}
  \put(-5,0){\line(1,0){24}}
  \put(-7,3){\line(1,0){28}}
\end{picture}                    & \\
\hline

$A_1$ & 1~~~~0 & \begin{picture}(5,15) \put(3,0){\circle*{6}} \put(5,5){\lab{6}} \end{picture} \\

$\widetilde{A}_1$ & 0~~~~1 &
\begin{picture}(5,15) \put(3,0){\circle{6}} \put(5,5){\lab{4}} \end{picture} \\

$A_1+\widetilde{A}_1$ & 2~~~~0 & 
\begin{picture}(5,15) 
\put(3,0){\circle*{6}} 
\put(5,5){\lab{2}} 
\put(23,0){\circle{6}}
\put(25,5){\lab{4}}
\end{picture} \\

$G_2$ & 2~~~~2 & 
\begin{picture}(10,17)
  \put(0,0){\circle{6}}
  \put(0,5){\lab{1}}
  \put(30,0){\circle*{6}}
  \put(32,5){\lab{2}}
  \put(1,-3){\line(1,0){28}}
  \put(3,0){\line(1,0){24}}
  \put(1,3){\line(1,0){28}}
\end{picture}                  \\
\hline

\end{longtable}

\begin{longtable}{|l|c|l|}
\caption{Nilpotent orbits in the Lie algebra of type $F_4$.}
\endfirsthead
\hline
\multicolumn{3}{|l|}{\small\slshape Nilpotent orbits in type $F_4$} \\ 
\hline
\endhead
\hline
\endfoot
\endlastfoot

\hline
label & diagram & representative \\

&  

\begin{picture}(20,7)
  \put(-20,0){\circle*{6}}
  \put(0,0){\circle*{6}}
  \put(20,0){\circle{6}}
  \put(40,0){\circle{6}}
  \put(-17,0){\line(1,0){14}}
  \put(2,2){\line(1,0){16}}
  \put(2,-2){\line(1,0){16}}
  \put(23,0){\line(1,0){14}}
\end{picture}                   

& \\
\hline

$A_1$ & 1~~~0~~~0~~~0 & \begin{picture}(5,15) \put(3,0){\circle*{6}} \put(5,5){\lab{24}} \end{picture} \\

$\widetilde{A}_1$ & 0~~~0~~~0~~~1 &
\begin{picture}(5,15) \put(3,0){\circle{6}} \put(5,5){\lab{21}} \end{picture} \\

$A_1+\widetilde{A}_1$ & 0~~~1~~~0~~~0 & 
\begin{picture}(5,15) 
\put(3,0){\circle{6}} 
\put(5,5){\lab{17}} 
\put(23,0){\circle*{6}}
\put(25,5){\lab{22}}
\end{picture} \\

$A_2$ & 2~~~0~~~0~~~0 & 
\begin{picture}(5,15) 
\put(3,0){\circle*{6}} 
\put(5,5){\lab{16}} 
\put(23,0){\circle*{6}}
\put(25,5){\lab{18}}
\put(5,0){\line(1,0){17}}
\end{picture} \\

$\widetilde{A}_2$ & 0~~~0~~~0~~~2 &
\begin{picture}(5,15) 
\put(3,0){\circle{6}} 
\put(5,5){\lab{11}} 
\put(23,0){\circle{6}}
\put(25,5){\lab{12}}
\put(6,0){\line(1,0){14}}
\end{picture} \\

$\widetilde{A}_1+A_2$ & 0~~~0~~~1~~~0 &
\begin{picture}(5,15)
\put(0,0){\circle{6}}
\put(0,5){\lab{14}} 
\put(18,0){\circle*{6}} 
\put(20,5){\lab{15}} 
\put(38,0){\circle*{6}}
\put(40,5){\lab{16}}
\put(21,0){\line(1,0){14}}
\end{picture}  \\
\hline
\end{longtable}

\begin{longtable}{|l|c|l|}
\caption{Nilpotent orbits in the Lie algebra of type $E_6$.}
\endfirsthead
\hline
\multicolumn{3}{|l|}{\small\slshape Nilpotent orbits in type $E_6$} \\ 
\hline
\endhead
\hline
\endfoot
\endlastfoot

\hline
label & diagram & representative \\

&  

\begin{picture}(40,25)
  \put(-20,0){\circle{6}}
  \put(0,0){\circle{6}}
  \put(20,0){\circle{6}}
  \put(40,0){\circle{6}}
  \put(60,0){\circle{6}}
  \put(20,20){\circle{6}}
  \put(-17,0){\line(1,0){14}}
  \put(3,0){\line(1,0){14}}
  \put(23,0){\line(1,0){14}}
  \put(43,0){\line(1,0){14}}
  \put(20,3){\line(0,1){14}}
\end{picture}                   

& \\
\hline

$A_1$ & $0~~~~0~~~~\overset{\text{\normalsize 1}}{0}~~~~0~~~~0$ & 
\begin{picture}(5,15) \put(3,0){\circle{6}} \put(5,5){\lab{36}} \end{picture} \\

$2A_1$ & $1~~~~0~~~~\overset{\text{\normalsize 0}}{0}~~~~0~~~~1$ & 
\begin{picture}(5,15) 
\put(3,0){\circle{6}} 
\put(5,5){\lab{30}} 
\put(23,0){\circle{6}}
\put(25,5){\lab{34}}
\end{picture} \\

$3A_1$ & $0~~~~0~~~~\overset{\text{\normalsize 0}}{1}~~~~0~~~~0$ & 
\begin{picture}(5,15) \put(3,0){\circle{6}} \put(5,5){\lab{24}} 
\put(23,0){\circle{6}}
\put(25,5){\lab{32}}
\put(43,0){\circle{6}}
\put(45,5){\lab{33}}
\end{picture} \\

$A_2$ & $0~~~~0~~~~\overset{\text{\normalsize 2}}{0}~~~~0~~~~0$ & 
\begin{picture}(5,15) 
\put(3,0){\circle{6}} 
\put(5,5){\lab{25}} 
\put(23,0){\circle{6}}
\put(25,5){\lab{26}}
\put(6,0){\line(1,0){14}}
\end{picture} \\

$A_2+A_1$ & $1~~~~0~~~~\overset{\text{\normalsize 1}}{0}~~~~0~~~~1$ & 
\begin{picture}(5,15) 
\put(3,0){\circle{6}} 
\put(5,5){\lab{25}} 
\put(23,0){\circle{6}}
\put(25,5){\lab{26}}
\put(6,0){\line(1,0){14}}
\put(43,0){\circle{6}}
\put(45,5){\lab{23}}
\end{picture} \\

$2A_2$ & $2~~~~0~~~~\overset{\text{\normalsize 0}}{0}~~~~0~~~~2$ & 
\begin{picture}(5,15) 
\put(3,0){\circle{6}} 
\put(5,5){\lab{17}} 
\put(23,0){\circle{6}}
\put(25,5){\lab{21}}
\put(6,0){\line(1,0){14}}
\put(43,0){\circle{6}}
\put(45,5){\lab{18}}
\put(63,0){\circle{6}}
\put(46,0){\line(1,0){14}}
\put(65,5){\lab{20}}
\end{picture} \\
\hline
\end{longtable}

\begin{longtable}{|l|c|l|}
\caption{Nilpotent orbits in the Lie algebra of type $E_7$.}
\endfirsthead
\hline
\multicolumn{3}{|l|}{\small Nilpotent orbits in $E_7$.} \\ 
\hline
\endhead
\hline
\endfoot
\endlastfoot

\hline
label & diagram & representative \\

&  

\begin{picture}(100,25)
  \put(0,0){\circle{6}}
  \put(20,0){\circle{6}}
  \put(40,0){\circle{6}}
  \put(60,0){\circle{6}}
  \put(80,0){\circle{6}}
  \put(100,0){\circle{6}}
  \put(40,20){\circle{6}}
  \put(3,0){\line(1,0){14}}
  \put(23,0){\line(1,0){14}}
  \put(43,0){\line(1,0){14}}
  \put(63,0){\line(1,0){14}}
  \put(83,0){\line(1,0){14}}
  \put(40,3){\line(0,1){14}}
\end{picture}                   

& \\
\hline

$A_1$ & $1~~~~0~~~~\overset{\text{\normalsize 0}}{0}~~~~0~~~~0~~~~0$ & 
\begin{picture}(5,15) \put(3,0){\circle{6}} \put(5,5){\lab{63}} \end{picture} \\

$2A_1$ & $0~~~~0~~~~\overset{\text{\normalsize 0}}{0}~~~~0~~~~1~~~~0$ & 
\begin{picture}(5,15) 
\put(3,0){\circle{6}} 
\put(5,5){\lab{57}} 
\put(23,0){\circle{6}}
\put(25,5){\lab{60}}
\end{picture} \\

$(3A_1)''$ & $0~~~~0~~~~\overset{\text{\normalsize 0}}{0}~~~~0~~~~0~~~~2$ & 
\begin{picture}(5,15) 
\put(3,0){\circle{6}} 
\put(5,5){\lab{47}} 
\put(23,0){\circle{6}}
\put(25,5){\lab{48}}
\put(43,0){\circle{6}}
\put(45,5){\lab{49}}
\end{picture} \\

$(3A_1)'$ & $0~~~~1~~~~\overset{\text{\normalsize 0}}{0}~~~~0~~~~0~~~~0$ & 
\begin{picture}(5,15) 
\put(3,0){\circle{6}} 
\put(5,5){\lab{42}} 
\put(23,0){\circle{6}}
\put(25,5){\lab{56}}
\put(43,0){\circle{6}}
\put(45,5){\lab{59}}
\end{picture} \\

$A_2$ & $2~~~~0~~~~\overset{\text{\normalsize 0}}{0}~~~~0~~~~0~~~~0$ & 
\begin{picture}(5,15) 
\put(3,0){\circle{6}} 
\put(5,5){\lab{44}} 
\put(23,0){\circle{6}}
\put(25,5){\lab{46}}
\put(6,0){\line(1,0){14}}
\end{picture} \\

$4A_1$ & $0~~~~0~~~~\overset{\text{\normalsize 1}}{0}~~~~0~~~~0~~~~1$ & 
\begin{picture}(5,15) 
\put(3,0){\circle{6}} 
\put(5,5){\lab{45}} 
\put(23,0){\circle{6}}
\put(25,5){\lab{47}}
\put(43,0){\circle{6}}
\put(45,5){\lab{52}}
\put(63,0){\circle{6}}
\put(65,5){\lab{53}}
\end{picture} \\
\hline
\end{longtable}

\begin{longtable}{|l|c|l|}
\caption{Nilpotent orbits in the Lie algebra of type $E_8$.}
\endfirsthead
\hline
\multicolumn{3}{|l|}{\small Nilpotent orbits in $E_8$.} \\ 
\hline
\endhead
\hline
\endfoot
\endlastfoot

\hline
label & diagram & representative \\

&  

\begin{picture}(120,25)
  \put(0,0){\circle{6}}
  \put(20,0){\circle{6}}
  \put(40,0){\circle{6}}
  \put(60,0){\circle{6}}
  \put(80,0){\circle{6}}
  \put(100,0){\circle{6}}
  \put(120,0){\circle{6}}
  \put(40,20){\circle{6}}
  \put(3,0){\line(1,0){14}}
  \put(23,0){\line(1,0){14}}
  \put(43,0){\line(1,0){14}}
  \put(63,0){\line(1,0){14}}
  \put(83,0){\line(1,0){14}}
  \put(40,3){\line(0,1){14}}
  \put(103,0){\line(1,0){14}}
\end{picture}                   

& \\
\hline

$A_1$ & $0~~~~0~~~~\overset{\text{\normalsize 0}}{0}~~~~0~~~~0~~~~0~~~~1$ & 
\begin{picture}(5,15) \put(3,0){\circle{6}} \put(5,5){\lab{120}} \end{picture} \\

$2A_1$ & $1~~~~0~~~~\overset{\text{\normalsize 0}}{0}~~~~0~~~~0~~~~0~~~~0$ & 
\begin{picture}(5,15) 
\put(3,0){\circle{6}} 
\put(5,5){\lab{113}} 
\put(23,0){\circle{6}}
\put(25,5){\lab{114}}
\end{picture} \\

$3A_1$ & $0~~~~0~~~~\overset{\text{\normalsize 0}}{0}~~~~0~~~~0~~~~1~~~~0$ & 
\begin{picture}(5,15) 
\put(3,0){\circle{6}} 
\put(5,5){\lab{104}} 
\put(23,0){\circle{6}}
\put(25,5){\lab{105}}
\put(43,0){\circle{6}}
\put(45,5){\lab{106}}
\end{picture} \\

$A_2$ & $0~~~~0~~~~\overset{\text{\normalsize 0}}{0}~~~~0~~~~0~~~~0~~~~2$ & 
\begin{picture}(5,15) 
\put(3,0){\circle{6}} 
\put(5,5){\lab{88}} 
\put(23,0){\circle{6}}
\put(25,5){\lab{90}}
\put(6,0){\line(1,0){14}}
\end{picture} \\

$4A_1$ & $0~~~~0~~~~\overset{\text{\normalsize 1}}{0}~~~~0~~~~0~~~~0~~~~0$ & 
\begin{picture}(50,15) 
\put(3,0){\circle{6}} 
\put(5,5){\lab{95}} 
\put(23,0){\circle{6}}
\put(25,5){\lab{97}}
\put(43,0){\circle{6}}
\put(45,5){\lab{98}}
\put(63,0){\circle{6}}
\put(65,5){\lab{103}}
\end{picture} \\

$A_2+A_1$ & $1~~~~0~~~~\overset{\text{\normalsize 0}}{0}~~~~0~~~~0~~~~0~~~~1$ & 
\begin{picture}(5,15) 
\put(3,0){\circle{6}} 
\put(5,5){\lab{88}} 
\put(23,0){\circle{6}}
\put(25,5){\lab{90}}
\put(6,0){\line(1,0){14}}
\put(43,0){\circle{6}}
\put(45,5){\lab{97}}
\end{picture} \\
\hline
\end{longtable}

\section{The exceptional root systems in {\sf GAP}}\label{sec:appB}

In this appendix we list the positive roots of the exceptional
root systems, in the order in which they appear in {\sf GAP}4. The 
tables have to be read from left to right, and from top to bottom.
So the first root is the one top left, the second root is the second one
on the first line, and the last root is the
one bottom right. For each root its coefficients with respect to
a basis of simple roots are given.\par
The roots for $F_4$ in Table \ref{tab:F4rt} may seem slightly strange.
This is due to the fact that in {\sf GAP}4 the positive roots are
ordered differently than usual. In this table the coefficients of 
each root with respect to the ``usual'' ordering of a basis of simple roots is given 
(i.e., as in \cite{bou4}). However, the roots are listed in the same order as
they are in {\sf GAP}4.

\begin{longtable}{|c|c|c|c|c|c|}
\caption{Positive roots in the root system of type $G_2$.}
\endfirsthead
\hline
\multicolumn{6}{|l|}{\small\slshape Positive roots in $G_2$} \\ 
\hline
\endhead
\hline
\endfoot
\endlastfoot

\hline

10 & 01 & 11 & 21 & 
31 & 32 \\
\hline
\end{longtable}

\begin{longtable}{|c|c|c|c|c|c|c|c|}
\caption{Positive roots in the root system of type $F_4$.}\label{tab:F4rt}
\endfirsthead
\hline
\multicolumn{8}{|l|}{\small\slshape Positive roots in $F_4$} \\ 
\hline
\endhead
\hline
\endfoot
\endlastfoot

\hline

0001 & 1000 & 0010 & 0100 & 
0011 & 1100 & 0110 & 0111\\
1110 & 0120 & 1111 & 0121 & 
1120 & 1121 & 0122 & 1220\\
1221 & 1122 & 1231 & 1222 & 
1232 & 1242 & 1342 & 2342\\

\hline
\end{longtable}

\begin{longtable}{|c|c|c|c|c|c|c|c|c|}
\caption{Positive roots in the root system of type $E_6$.}
\endfirsthead
\hline
\multicolumn{9}{|l|}{\small\slshape Positive roots in $E_6$} \\ 
\hline
\endhead
\hline
\endfoot
\endlastfoot

\hline

$10\overset{\text{\normalsize 0}}{0}00$ & $00\overset{\text{\normalsize 1}}{0}00$ & 
$01\overset{\text{\normalsize 0}}{0}00$ & $00\overset{\text{\normalsize 0}}{1}00$ & 
$00\overset{\text{\normalsize 0}}{0}10$ & $00\overset{\text{\normalsize 0}}{0}01$ & 
$11\overset{\text{\normalsize 0}}{0}00$ & $00\overset{\text{\normalsize 1}}{1}00$ & 
$01\overset{\text{\normalsize 0}}{1}00$\\$00\overset{\text{\normalsize 0}}{1}10$ & 
$00\overset{\text{\normalsize 0}}{0}11$ & $11\overset{\text{\normalsize 0}}{1}00$ & 
$01\overset{\text{\normalsize 1}}{1}00$ & $00\overset{\text{\normalsize 1}}{1}10$ & 
$01\overset{\text{\normalsize 0}}{1}10$ & $00\overset{\text{\normalsize 0}}{1}11$ & 
$11\overset{\text{\normalsize 1}}{1}00$ & $11\overset{\text{\normalsize 0}}{1}10$\\
$01\overset{\text{\normalsize 1}}{1}10$ & $00\overset{\text{\normalsize 1}}{1}11$ & 
$01\overset{\text{\normalsize 0}}{1}11$ & $11\overset{\text{\normalsize 1}}{1}10$ & 
$11\overset{\text{\normalsize 0}}{1}11$ & $01\overset{\text{\normalsize 1}}{2}10$ & 
$01\overset{\text{\normalsize 1}}{1}11$ & $11\overset{\text{\normalsize 1}}{2}10$ & 
$11\overset{\text{\normalsize 1}}{1}11$\\$01\overset{\text{\normalsize 1}}{2}11$ & 
$12\overset{\text{\normalsize 1}}{2}10$ & $11\overset{\text{\normalsize 1}}{2}11$ & 
$01\overset{\text{\normalsize 1}}{2}21$ & $12\overset{\text{\normalsize 1}}{2}11$ & 
$11\overset{\text{\normalsize 1}}{2}21$ & $12\overset{\text{\normalsize 1}}{2}21$ & 
$12\overset{\text{\normalsize 1}}{3}21$ & $12\overset{\text{\normalsize 2}}{3}21$\\

\hline

\end{longtable}

\begin{longtable}{|c|c|c|c|c|c|c|}
\caption{Positive roots in the root system of type $E_7$.}
\endfirsthead
\hline
\multicolumn{7}{|l|}{\small\slshape Positive roots in $E_7$} \\ 
\hline
\endhead
\hline
\endfoot
\endlastfoot

\hline

$10\overset{\text{\normalsize 0}}{0}000$ & $00\overset{\text{\normalsize 1}}{0}000$ & 
$01\overset{\text{\normalsize 0}}{0}000$ & $00\overset{\text{\normalsize 0}}{1}000$ & 
$00\overset{\text{\normalsize 0}}{0}100$ & $00\overset{\text{\normalsize 0}}{0}010$ & 
$00\overset{\text{\normalsize 0}}{0}001$\\$11\overset{\text{\normalsize 0}}{0}000$ & 
$00\overset{\text{\normalsize 1}}{1}000$ & $01\overset{\text{\normalsize 0}}{1}000$ & 
$00\overset{\text{\normalsize 0}}{1}100$ & $00\overset{\text{\normalsize 0}}{0}110$ & 
$00\overset{\text{\normalsize 0}}{0}011$ & $11\overset{\text{\normalsize 0}}{1}000$\\
$01\overset{\text{\normalsize 1}}{1}000$ & $00\overset{\text{\normalsize 1}}{1}100$ & 
$01\overset{\text{\normalsize 0}}{1}100$ & $00\overset{\text{\normalsize 0}}{1}110$ & 
$00\overset{\text{\normalsize 0}}{0}111$ & $11\overset{\text{\normalsize 1}}{1}000$ & 
$11\overset{\text{\normalsize 0}}{1}100$\\$01\overset{\text{\normalsize 1}}{1}100$ & 
$00\overset{\text{\normalsize 1}}{1}110$ & $01\overset{\text{\normalsize 0}}{1}110$ & 
$00\overset{\text{\normalsize 0}}{1}111$ & $11\overset{\text{\normalsize 1}}{1}100$ & 
$11\overset{\text{\normalsize 0}}{1}110$ & $01\overset{\text{\normalsize 1}}{2}100$\\
$01\overset{\text{\normalsize 1}}{1}110$ & $00\overset{\text{\normalsize 1}}{1}111$ & 
$01\overset{\text{\normalsize 0}}{1}111$ & $11\overset{\text{\normalsize 1}}{2}100$ & 
$11\overset{\text{\normalsize 1}}{1}110$ & $11\overset{\text{\normalsize 0}}{1}111$ & 
$01\overset{\text{\normalsize 1}}{2}110$\\$01\overset{\text{\normalsize 1}}{1}111$ & 
$12\overset{\text{\normalsize 1}}{2}100$ & $11\overset{\text{\normalsize 1}}{2}110$ & 
$11\overset{\text{\normalsize 1}}{1}111$ & $01\overset{\text{\normalsize 1}}{2}210$ & 
$01\overset{\text{\normalsize 1}}{2}111$ & $12\overset{\text{\normalsize 1}}{2}110$\\
$11\overset{\text{\normalsize 1}}{2}210$ & $11\overset{\text{\normalsize 1}}{2}111$ & 
$01\overset{\text{\normalsize 1}}{2}211$ & $12\overset{\text{\normalsize 1}}{2}210$ & 
$12\overset{\text{\normalsize 1}}{2}111$ & $11\overset{\text{\normalsize 1}}{2}211$ & 
$01\overset{\text{\normalsize 1}}{2}221$\\$12\overset{\text{\normalsize 1}}{3}210$ & 
$12\overset{\text{\normalsize 1}}{2}211$ & $11\overset{\text{\normalsize 1}}{2}221$ & 
$12\overset{\text{\normalsize 2}}{3}210$ & $12\overset{\text{\normalsize 1}}{3}211$ & 
$12\overset{\text{\normalsize 1}}{2}221$ & $12\overset{\text{\normalsize 2}}{3}211$\\
$12\overset{\text{\normalsize 1}}{3}221$ & $12\overset{\text{\normalsize 2}}{3}221$ & 
$12\overset{\text{\normalsize 1}}{3}321$ & $12\overset{\text{\normalsize 2}}{3}321$ & 
$12\overset{\text{\normalsize 2}}{4}321$ & $13\overset{\text{\normalsize 2}}{4}321$ & 
$23\overset{\text{\normalsize 2}}{4}321$\\

\hline
\end{longtable}

\begin{longtable}{|c|c|c|c|c|c|}
\caption{Positive roots in the root system of type $E_8$.}
\endfirsthead
\hline
\multicolumn{6}{|l|}{\small\slshape Positive roots in $E_8$} \\ 
\hline
\endhead
\hline
\endfoot
\endlastfoot

\hline

$10\overset{\text{\normalsize 0}}{0}0000$ & $00\overset{\text{\normalsize 1}}{0}0000$ & 
$01\overset{\text{\normalsize 0}}{0}0000$ & $00\overset{\text{\normalsize 0}}{1}0000$ & 
$00\overset{\text{\normalsize 0}}{0}1000$ & $00\overset{\text{\normalsize 0}}{0}0100$\\
$00\overset{\text{\normalsize 0}}{0}0010$ & $00\overset{\text{\normalsize 0}}{0}0001$ & 
$11\overset{\text{\normalsize 0}}{0}0000$ & $00\overset{\text{\normalsize 1}}{1}0000$ & 
$01\overset{\text{\normalsize 0}}{1}0000$ & $00\overset{\text{\normalsize 0}}{1}1000$\\
$00\overset{\text{\normalsize 0}}{0}1100$ & $00\overset{\text{\normalsize 0}}{0}0110$ & 
$00\overset{\text{\normalsize 0}}{0}0011$ & $11\overset{\text{\normalsize 0}}{1}0000$ & 
$01\overset{\text{\normalsize 1}}{1}0000$ & $00\overset{\text{\normalsize 1}}{1}1000$\\
$01\overset{\text{\normalsize 0}}{1}1000$ & $00\overset{\text{\normalsize 0}}{1}1100$ & 
$00\overset{\text{\normalsize 0}}{0}1110$ & $00\overset{\text{\normalsize 0}}{0}0111$ & 
$11\overset{\text{\normalsize 1}}{1}0000$ & $11\overset{\text{\normalsize 0}}{1}1000$\\
$01\overset{\text{\normalsize 1}}{1}1000$ & $00\overset{\text{\normalsize 1}}{1}1100$ & 
$01\overset{\text{\normalsize 0}}{1}1100$ & $00\overset{\text{\normalsize 0}}{1}1110$ & 
$00\overset{\text{\normalsize 0}}{0}1111$ & $11\overset{\text{\normalsize 1}}{1}1000$\\
$11\overset{\text{\normalsize 0}}{1}1100$ & $01\overset{\text{\normalsize 1}}{2}1000$ & 
$01\overset{\text{\normalsize 1}}{1}1100$ & $00\overset{\text{\normalsize 1}}{1}1110$ & 
$01\overset{\text{\normalsize 0}}{1}1110$ & $00\overset{\text{\normalsize 0}}{1}1111$\\
$11\overset{\text{\normalsize 1}}{2}1000$ & $11\overset{\text{\normalsize 1}}{1}1100$ & 
$11\overset{\text{\normalsize 0}}{1}1110$ & $01\overset{\text{\normalsize 1}}{2}1100$ & 
$01\overset{\text{\normalsize 1}}{1}1110$ & $00\overset{\text{\normalsize 1}}{1}1111$\\
$01\overset{\text{\normalsize 0}}{1}1111$ & $12\overset{\text{\normalsize 1}}{2}1000$ & 
$11\overset{\text{\normalsize 1}}{2}1100$ & $11\overset{\text{\normalsize 1}}{1}1110$ & 
$11\overset{\text{\normalsize 0}}{1}1111$ & $01\overset{\text{\normalsize 1}}{2}2100$\\
$01\overset{\text{\normalsize 1}}{2}1110$ & $01\overset{\text{\normalsize 1}}{1}1111$ & 
$12\overset{\text{\normalsize 1}}{2}1100$ & $11\overset{\text{\normalsize 1}}{2}2100$ & 
$11\overset{\text{\normalsize 1}}{2}1110$ & $11\overset{\text{\normalsize 1}}{1}1111$\\
$01\overset{\text{\normalsize 1}}{2}2110$ & $01\overset{\text{\normalsize 1}}{2}1111$ & 
$12\overset{\text{\normalsize 1}}{2}2100$ & $12\overset{\text{\normalsize 1}}{2}1110$ & 
$11\overset{\text{\normalsize 1}}{2}2110$ & $11\overset{\text{\normalsize 1}}{2}1111$\\
$01\overset{\text{\normalsize 1}}{2}2210$ & $01\overset{\text{\normalsize 1}}{2}2111$ & 
$12\overset{\text{\normalsize 1}}{3}2100$ & $12\overset{\text{\normalsize 1}}{2}2110$ & 
$12\overset{\text{\normalsize 1}}{2}1111$ & $11\overset{\text{\normalsize 1}}{2}2210$\\
$11\overset{\text{\normalsize 1}}{2}2111$ & $01\overset{\text{\normalsize 1}}{2}2211$ & 
$12\overset{\text{\normalsize 2}}{3}2100$ & $12\overset{\text{\normalsize 1}}{3}2110$ & 
$12\overset{\text{\normalsize 1}}{2}2210$ & $12\overset{\text{\normalsize 1}}{2}2111$\\
$11\overset{\text{\normalsize 1}}{2}2211$ & $01\overset{\text{\normalsize 1}}{2}2221$ & 
$12\overset{\text{\normalsize 2}}{3}2110$ & $12\overset{\text{\normalsize 1}}{3}2210$ & 
$12\overset{\text{\normalsize 1}}{3}2111$ & $12\overset{\text{\normalsize 1}}{2}2211$\\
$11\overset{\text{\normalsize 1}}{2}2221$ & $12\overset{\text{\normalsize 2}}{3}2210$ & 
$12\overset{\text{\normalsize 2}}{3}2111$ & $12\overset{\text{\normalsize 1}}{3}3210$ & 
$12\overset{\text{\normalsize 1}}{3}2211$ & $12\overset{\text{\normalsize 1}}{2}2221$\\
$12\overset{\text{\normalsize 2}}{3}3210$ & $12\overset{\text{\normalsize 2}}{3}2211$ & 
$12\overset{\text{\normalsize 1}}{3}3211$ & $12\overset{\text{\normalsize 1}}{3}2221$ & 
$12\overset{\text{\normalsize 2}}{4}3210$ & $12\overset{\text{\normalsize 2}}{3}3211$\\
$12\overset{\text{\normalsize 2}}{3}2221$ & $12\overset{\text{\normalsize 1}}{3}3221$ & 
$13\overset{\text{\normalsize 2}}{4}3210$ & $12\overset{\text{\normalsize 2}}{4}3211$ & 
$12\overset{\text{\normalsize 2}}{3}3221$ & $12\overset{\text{\normalsize 1}}{3}3321$\\
$23\overset{\text{\normalsize 2}}{4}3210$ & $13\overset{\text{\normalsize 2}}{4}3211$ & 
$12\overset{\text{\normalsize 2}}{4}3221$ & $12\overset{\text{\normalsize 2}}{3}3321$ & 
$23\overset{\text{\normalsize 2}}{4}3211$ & $13\overset{\text{\normalsize 2}}{4}3221$\\
$12\overset{\text{\normalsize 2}}{4}3321$ & $23\overset{\text{\normalsize 2}}{4}3221$ & 
$13\overset{\text{\normalsize 2}}{4}3321$ & $12\overset{\text{\normalsize 2}}{4}4321$ & 
$23\overset{\text{\normalsize 2}}{4}3321$ & $13\overset{\text{\normalsize 2}}{4}4321$\\
$23\overset{\text{\normalsize 2}}{4}4321$ & $13\overset{\text{\normalsize 2}}{5}4321$ & 
$23\overset{\text{\normalsize 2}}{5}4321$ & $13\overset{\text{\normalsize 3}}{5}4321$ & 
$23\overset{\text{\normalsize 3}}{5}4321$ & $24\overset{\text{\normalsize 2}}{5}4321$\\
$24\overset{\text{\normalsize 3}}{5}4321$ & $24\overset{\text{\normalsize 3}}{6}4321$ & 
$24\overset{\text{\normalsize 3}}{6}5321$ & $24\overset{\text{\normalsize 3}}{6}5421$ & 
$24\overset{\text{\normalsize 3}}{6}5431$ & $24\overset{\text{\normalsize 3}}{6}5432$\\

\hline

\end{longtable}


\def\cprime{$'$} \def\cprime{$'$} \def\cprime{$'$} \def\cprime{$'$}

\end{document}